\documentclass[10pt, reqno]{amsart}

\usepackage{graphicx}

\usepackage{amssymb, amsfonts, verbatim}

\usepackage{epsfig, amsmath, cancel}

\usepackage[mathscr]{eucal}

\newtheorem{lemma}{Lemma}

\newtheorem{prop}{Proposition}

\newtheorem{theorem}{Theorem}

\newcommand{\Q}{\hspace{.065in}} 

\newcommand{\diag}[1]{\ensuremath{\operatorname{diag}}\left( #1 \right)}

\newcommand{\delmu}{\partial_{\mu}}

\newcommand{\delnu}{\partial_{\nu}}

\newcommand{\dela}{\partial_{\alpha}}

\newcommand{\delb}{\partial_{\beta}}

\newcommand{\del}[1]{\partial_{#1}}

\newcommand{\Del}[1]{\partial^{#1}}

\newcommand{\Lonet}[2]{ \Big\| #2 (#1,\cdot) \Big\|_{L^1} }

\newcommand{\ltwot}[2]{ \| #2 (#1,\cdot) \|_{L^2} }

\newcommand{\linftyt}[2]{ \| #2 (#1,\cdot) \|_{L^{\infty}} }

\newcommand{\ga}{\alpha}

\newcommand{\gb}{\beta}

\newcommand{\gd}{\delta}

\renewcommand{\ge}{\varepsilon}

\renewcommand{\gg}{\gamma} 

\newcommand{\gh}{\eta}

\newcommand{\gm}{\mu}

\newcommand{\gn}{\nu}

\newcommand{\gt}{\tau}

\newcommand{\gy}{\psi}

\newcommand{\gO}{\Omega}

\newcommand{\scL}{\mathscr{L}}

\newcommand{\bbR}{\mathbb{R}}

\newcommand{\calO}{\mathcal{O}}

\begin{document}

\title[Timelike Minimal Surfaces]{Timelike Minimal Submanifolds of General Co-dimension in Minkowski SpaceTime}
\author{Paul Allen, Lars Andersson, James Isenberg}
\address{ Department of Mathematics\\ 
	University of Oregon\\
	Eugene, OR 97403, USA}
\email{pallen1@math.uoregon.edu}
\address{AEI\\ 
	Max Planck Institute for Gravitational Physics\\
	D-14476 Potsdam, Germany \and 
Department of Mathematics\\
University of Miami\\
Coral Gables, FL 33124\\
USA}
\email{andersson@aei.mpg.de}
\address{ Department of Mathematics\\ 
	University of Oregon\\
	Eugene, OR 97403, USA}
\email{jim@newton.uoregon.edu}
\date{ 1 December 2005}

\begin{abstract}
We consider the timelike minimal surface problem in Minkowski spacetimes and
show local and global existence of such surfaces having arbitrary dimension
$\geq 2$ and arbitrary 
co-dimension, provided they are initially close to a flat plane.
\end{abstract}

\maketitle

\numberwithin{equation}{section} 


\section{Introduction}

In this work we consider timelike minimal submanifolds of dimension $1+n$,
$n\geq 2$, of 
Minkowski spacetimes of dimension $1+n+q$, $q \geq 1$. A submanifold is called minimal if
it is stationary with respect to variations of the induced area, which
thus provides an action for the system. Timelike minimal submanifolds may be
viewed as simple but nontrivial examples of D-branes, which play an important
role in string-theory, and the system under consideration here thus has
natural generalizations motivated by string theory. 
In this work we prove a small data, global existence result for timelike
minimal submanifolds of arbitrary codimension $q$. The solutions which are
constructed are close to flat timelike planes. 

The Euler-Lagrange equations arising from variation of the area 
form a quasilinear system of PDE's, 
which under suitable conditions on the data is hyperbolic. This system is
closely related to the scalar quasilinear hyperbolic PDE governing 
timelike minimal
hypersurfaces. The small data, global existence problem for timelike minimal
hypersurfaces has been considered by Brendle \cite{B_Minimal} and Lindblad \cite{L_MinSurf}. The
work of Lindblad 
makes use of the null structure of the system, and our
approach is closely related to the work in \cite{L_MinSurf}.

Consider an embedding of $\bbR^{1+n}$ into Minkowski spacetime $\bbR^{1+n+q}$
given by the graph of a map $f: \bbR^{1+n} \to \bbR^q$. Let greek indices
$\alpha, \beta, \dots$ take values in $0,1,\dots,n$ and let uppercase latin
indices $I,J,\dots$ take values in $1,\dots,q$. We introduce cartesian
coordinates $x^\alpha$ on $\bbR^{1+n}$ and $x^I$ on $\bbR^q$. 
The induced metric $\bbR^{1+n}$ is
	\begin{equation}
	 h_{\ga\gb} = \gh_{\ga\gb} + f^I_{\ga}f^J_{\gb}\gd_{IJ} ,
	\end{equation}
where $f^I = x^I\circ f$, $f^I_{\ga} = \dela f^I$ and $\gh = \diag{-1,1\dots,1}$ is the Minkowski
metric.  Varying the action 
$$
\mathcal S = \int
\sqrt{- \det h_{\mu\nu}} \, d^{1+n} x, 
$$
yields the
Euler-Lagrange equations
	\begin{equation} \label{ELEqn}
	0 = \delmu \left[ \sqrt{-\det{h}}\, h^{\gm\gn} f^I_{\gn}\right]
            \quad I=1,\dots,q ,
	\end{equation}
which we consider for small data
	\begin{equation}\label{E:Data}
	f^I(0,\cdot) = \ge g^I \quad \del{t}f^I(0,\cdot) = \ge k^I,
	\end{equation}	
with $g^I,k^I$ smooth and decaying sufficiently fast for large $|x|$; for
simplicity we restrict to $\ge \leq 1$.  Here $h^{\mu \nu}$ is the inverse of
$h_{\alpha \beta}$.\\
	
For future use, we note that equation (\ref{ELEqn}) can be written 
in divergence form 
	\begin{equation}\label{E:DivForm}
	\Box f^I = \delmu \left[ F^{\gm\gn} f^I_{\gn}\right] , 
	\end{equation}
where $\Box = \gh^{\gm\gn}\delmu\delnu$ is the Minkowski wave operator and
	$
	F^{\gm\gn}(\del{}f) = \gh^{\gm\gn} -  \sqrt{-\det{h}} \, h^{\gm\gn}, 
	$
as well as in the form 
	\begin{equation}\label{E:SymSystem}
	H^{\gm\gn}_{\Q\Q JL} (\del{}f) \delmu\delnu f^J = 0, \qquad I=1,\dots,q 
	\end{equation}
where
	\begin{equation}\label{E:Hdef}
	H^{\gm\gn}_{\Q\Q JL} 
	= \sqrt{-\det{h}}\left[ \gd_{JL} h^{\gm\gn} 
		- \gd_{IJ}\gd_{KL} \left( h^{\gm\gn} h^{\ga\gb}f_{\ga}^Kf_{\gb}^I 
			+ h^{\gm\ga}h^{\gn\gb} f_{\ga}^I f_{\gb}^K
			+ h^{\gm\ga}h^{\gn\gb} f_{\ga}^K f_{\gb}^I
		 \right) \right] .
	\end{equation}
We raise and lower Greek (intrinsic) indices using $h_{\mu \nu}$ and its
inverse, while Latin (extrinsic) indices are raised and lowered using the
identity $\delta_{I J}$ and its inverse.  From equation 
\eqref{E:Hdef} it follows that $H^{\gm\gn}_{\Q\Q JL}$ has the symmetries 
\begin{equation}\label{E:symm} 
H^{\gm\gn}_{\Q\Q JL} = H^{\gm\gn}_{\Q\Q LJ} = H^{\gn\gm}_{\Q\Q JL}.
\end{equation}
Due to the symmetries, an energy estimate and local well posedness 
holds for the system (\ref{E:SymSystem}).

The local existence argument follows \cite{H_Book} (see also
\cite{So_Lectures}) and uses an energy inequality which takes advantage of
symmetries in the system; see \cite{Keel} for a treatment of well-posedness
and lifespan results for symmetric systems.  The global existence result also
requires estimates applicable to divergence equations (see \cite{L_Lifespan})
and 
an
$L^{\infty} - L^1$ estimate.  Furthermore, the global existence result in $n=2$
dimension exploits the fact that the equation satisfies the null
condition; see Section \ref{S:DimTwo} for details.

\section{Local Existence}\label{S:MinkLocalExist}

The local well-posedness for systems of the form (\ref{E:SymSystem}) is well
established. Therefore we do not give a complete proof here but rather provide a
simple proof of the basic energy estimate.
The proof of local well-posedness
then follows along the same lines as the proofs given in \cite{H_Book} or 
\cite{So_Lectures}. Note that the energy estimate we state here also plays a key role in the global existence results discussed in the next section.

Write $|\del{}f|^2 = \gh^{\gm\gn}\gd_{IJ}f^I_{\gm}f^J_{\gn}$; the function
space norms used below are defined in terms of this and analogous
expressions. 
Using the identity 
\begin{multline*}
	(f^I_0)({H}^{\gm\gn}_{\Q\Q IJ} f^J_{\gm\gn}) = \delmu \left[
	{H}^{\gm\gn}_{\Q\Q IJ} f^I_{\gn} f^J_0 - \frac12
	{H}^{\ga\gb}_{\Q\Q IJ} f^I_{\ga}f^J_{\gb} \gd^\gm_0 \right] 
\\ 
       -\delmu\left[ H^{\gm\gn}_{\Q\Q IJ}\right] f^I_{\gn} f^J_0 +
	\frac12 \del{0}\left[ H^{\gm\gn}_{\Q\Q IJ}\right] f^I_{\gm}
	f^J_{\gn},
	\end{multline*}
a standard argument yields the following energy estimate. 
\begin{lemma}\label{L:SymEnergy}
Let $f\in C^2([0,T)\times \bbR^n; \bbR^q)$ 
for $T > 0$ 
and assume that $H^{\gm\gn}_{\Q\Q IJ}$ has the symmetries 
(\ref{E:symm}). 
Further, assume that $H^{\gm\gn}_{\Q\Q IJ}$ 
satisfies
	\begin{equation}\label{E:coerce} 
	\sum \lvert H^{\gm\gn}_{\Q\Q IJ} - \gh^{\gm\gn}\gd_{IJ}\rvert
	< \frac12 \qquad \text{ on } [0,T].
	\end{equation} 
Then for $t\in [0,T)$ we have
	\begin{multline}\label{E:SymEnergy}
	 \ltwot{t}{\del{}f} \leq 2 \left( \ltwot{0}{\del{}f} + \int_0^t
	 \sum_{I} \ltwot{\gt}{H^{\gm\gn}_{\Q\Q IJ}\delmu\delnu
	 f^J}d\gt\right)\\ \times \exp{\left(\int_0^t 2
	 \linftyt{s}{\del{} H} ds\right)}.
	\end{multline}
\end{lemma}
With this energy estimate, the proof of local well-posedness can now be
completed by an iteration procedure following exactly the outline in
\cite{So_Lectures} or \cite{H_Book}. 
The approach in \cite{H_Book} makes use of some
extra structural assumptions which are easily removed, and gives local
well-posedness in Sobolev spaces $H^{s}$, for integer $s > n/2+2$. The
argument in Sogge gives the result for $s > n+3$. Since we are concerned
with small data, global existence here, the exact regularity needed for the
local well-posedness is not important. We can now state the following result. 

\begin{theorem} \label{T:cauchy}
Let $s > n/2 +2$ and consider equation 
\begin{subequations}\label{E:cauchy} 
\begin{align}
H^{\gm\gn}_{\Q\Q JL} (x,f, \partial f) 
\delmu\delnu  f^J &= G_L(x,f,\partial f), \qquad I=1,\dots,q \\
\intertext{with initial data}
f(0,\cdot) = g, \qquad \partial_0 f(0,\cdot) &= k
\end{align}
\end{subequations}
Here $H^{\gm\gn}_{\Q\Q JL}$ and $G_L$ are assumed to be smooth functions of
their arguments and  $H^{\gm\gn}_{\Q\Q JL}$ is assumed to satisfy 
the symmetries
\eqref{E:symm}.  
Suppose the initial data $(g ,k)$ is such that equation
\eqref{E:coerce} is valid for $H^{\gm\gn}_{\Q\Q JL}$ evaluated on $(g, k)$. 
Then there is a $T > 0$, which depends only on the norm of
$(g, k)$ in $H^s\times H^{s-1}$,
and a function $f \in C^2([0,T]\times \bbR^n;\bbR^q)$ 
which solves (\ref{E:cauchy}), with $|\partial^\alpha
f|$ bounded for $|\alpha| \leq 2$. The maximal time of existence is bounded
from below by the
supremum of all $T$ such that (\ref{E:cauchy}) has a $C^2$ solution such that
for $0 \leq t \leq T$, equation (\ref{E:coerce}) is valid and 
$\partial^\alpha f$ is bounded for $|\alpha| \leq 2$.  
\end{theorem}

\section{Global existence in dimensions $n\geq 3$}\label{S:GlobeExist}

Global existence follows from a procedure similar to that discussed in
\cite{L_MinSurf}.  The estimates needed rely on a collection of weighted norms
defined using the set of Lorentz vector fields $Z_{\gm}$, which include the
generators of the Lorentz group, along with the generator of dilations:
\begin{equation*}
	\{ \delmu, \gO_{ab}:=x^b\del{a} - x^a\del{b},
	\gO_{0a}:=t\del{a}+x^a\del{t}, L:=t\del{t}+r\del{r}\}.
\end{equation*}
These form a Lie algebra, which satisfies the following commutation relations
	\begin{equation*}
	[Z_{\gm}, \delnu ] = \sum a_{\gm\gn}^{\Q\Q \ga}\dela, \quad
		a_{\gm\gn}^{\Q\Q \ga} = 0,\pm 1,\quad \text{and}\quad [Z_\gm,
		\Box] = \begin{cases} -2\Box & \text{ if } Z_{\gm}=L,\\ 0 &
		\text{ otherwise.}\end{cases}
	\end{equation*}
It follows that  the equation $\Box \gy =0$ is preserved by these operators.
We now define the following norms in terms of products $Z^{\ga}$ (for
multi-index $\ga$) of the Lorentz vector fields applied to $f$:
\begin{subequations}\label{E:Norms} 
\begin{align}
M_1(t) &\equiv \sum_{|\ga|\leq m} \ltwot{t}{\del{}Z^{\ga} f} ,\\ M_2(t)
&\equiv \sum_{|\ga|\leq m} \ltwot{t}{Z^{\ga}f} , \\ N_1(t) &\equiv
\sum_{|\ga|\leq l} \linftyt{t}{\del{}Z^{\ga}f} ,\\ N_2(t) &\equiv
\sum_{|\ga|\leq l+1} \linftyt{t}{Z^{\ga}f} ,
\end{align} 
\end{subequations} 
Here $m$ is an integer such that $m > 2n+1$ and $l = (m+1)/2$.
Since $\delmu \in \{ Z_{\gn}\}$, these norms control the $L^2$ and
$L^{\infty}$ Sobolev norms (of order $k$ and $l$, respectively) of $f$ and
$\del{}f$. It follows that controlling these norms is sufficient to overcome
the obstruction to local existence discussed at the end of the previous
section.  In particular, if we control these norms, then $f$ is bounded in
$C^2$. Furthermore, we use the norms $N_1, N_2$ to show stability in the
sense that they decay as $t$ grows.  Note that any estimate for $N_2$ implies
an estimate for $N_1$ as well.\\


The following three propositions also play a role in our proof of global
existence. (See \cite{L_Lifespan}, and also 
\cite{L_MinSurf},\cite{Hor_Linfty},\cite{Kl_Linfty}, for some of the  proofs
of these results.)

\begin{prop}\label{P:DivIeq}
If $g$ is a solution to
	\begin{equation*}
	\begin{cases}
	 \Box g = \delmu F^{\gm} \\ g(0,\cdot) = \ge g, \quad
	 \del{0}g(0,\cdot) =\ge k
	\end{cases} 
	\end{equation*}
then
	\begin{equation*}
	 \ltwot{t}{g} \leq C_{\text{data}}\ge \, m(t) + \sum_{\gm} \int_0^t
	 \ltwot{\gt}{F^{\gm}} \,d\gt,
	\end{equation*}
where
	\begin{equation*}
	m(t) = \begin{cases} \log{(2+t)} & \text{ if } n=2,\\ 1 & \text{
 		 otherwise, } \end{cases}
	\end{equation*}
and $C_{\text{data}}$ depends on $g, k$, and  $F^0(0,\cdot)$.
\end{prop}

\begin{prop}\label{P:Linfty}
The solution to
	\begin{equation*}
	 \begin{cases}
	 \Box g = G \\ g(0,\cdot) = \ge g_0, \quad \del{0}g(0,\cdot) =\ge g_1
	\end{cases} 
	\end{equation*}
satisfies
	\begin{equation*}
	 \lvert g(t,x) \rvert \leq \frac{C}{(1+t+|x|)^{(n-1)/2}} \left(
	 	C_{\text{data}} \ge  + \int_0^t \sum_{|\ga|\leq n-1}
	 	\Lonet{s}{\frac{Z^{\ga}G}{(1+s+|\cdot|)^{(n-1)/2}}}\,ds
	 	\right).
	\end{equation*}\\
\end{prop}


To make use of the above propositions in controlling the norms
\eqref{E:Norms}, we apply $Z^{\ga}$ to both sides of the equation
\eqref{E:SymSystem} for $f^I$, and obtain 
	\begin{equation}\label{E:SymZ}
	{H}^{\gm\gn}_{\Q\Q IJ} \delmu\delnu (Z^{\ga}f^J) =  \sum_{k\geq 3,
		{\sum |\ga_i|\leq |\ga|+1}  }  H_{I, I_1\cdots I_k,
		\gg_1\cdots \gg_k,\ga_1\cdots \ga_k}\,
		(\del{\gg_1}Z^{\ga_1}f^{I_1})\cdots
		(\del{\gg_k}Z^{\ga_k}f^{I_k}).
	\end{equation}
Since $h^{\gm\gn} = \gh^{\gm\gn} + \calO(|\del{}f|^2)$, this may also be
written as
	\begin{equation}\label{E:BoxZ}
	 \Box \left( Z^{\ga} f^I\right)  =\sum_{k\geq 3, {\sum |\ga_i|\leq
		|\ga|+1} }  \hat{H}^I_{I_1\cdots I_k, \gg_1\cdots
		\gg_k,\ga_1\cdots \ga_k}\,
		(\del{\gg_1}Z^{\ga_1}f^{I_1})\cdots
		(\del{\gg_k}Z^{\ga_k}f^{I_k}).
\end{equation}
with some modified coefficient functions $\hat{H}$, satisfying $\hat{H} =
\calO(|\partial f|^2)$. Note that  at most one of the $\alpha_i$ can satisfy
$|\alpha_i| > (|\alpha|+1)/2$.  The global existence proof also depends on
the form of the divergence equation \eqref{E:DivForm}.  In particular, we
note that $ \sqrt{\det{[h]}}\,h^{\gm\gn} = \gh^{\gm\gn} + \calO(|\del{}f|^2)
$ as $|\del{}f|\to 0$, and thus $F^{\gm\gn} :=  \gh^{\gm\gn} -  \sqrt{-\det{h}} \, h^{\gm\gn}=\calO(|\del{}f|^2)$.  Hence
applying $Z^{\ga}$ to \eqref{E:DivForm} we obtain
	\begin{equation}\label{E:DivZ}
	\Box \left(Z^{\ga} f^I\right)  =\delmu \left[ \sum_{k\geq 3, \sum
		|\ga_i|\leq |\ga| } F^{\gm,I}_{I_1\cdots I_k,
		\gg_1\cdots\gg_k, \ga_1\cdots \ga_k}\,
		(\del{\gg_1}Z^{\ga_1}f^{I_1}) \cdots(\del{\gg_k}
		Z^{\ga_k}f^{I_k}) \right] .
	\end{equation}
where again  at most one of the $\alpha_i$ can satisfy $|\alpha_i| >
|\alpha|/2$.

We are now prepared to
show global existence in dimensions $n\geq 3$. The proof uses a continuous
induction, or bootstrap argument.
%


To set up the bootstrap argument, we assume that there is a constant $K$ so
that on $[0,T)$ we have the following estimates for the norms defined in
\eqref{E:Norms}:
\begin{subequations}\label{E:NormBd}
\begin{align}
M_1(t)  &\leq K\ge, \\
M_2(t)  &\leq K\ge, \\
N_1(t)  &\leq \frac{K\ge}{(1+t)^{\frac{n-1}{2}}} ,\\
N_2(t)  &\leq \frac{K\ge}{(1+t)^{\frac{n-1}{2}}}.
\end{align} 
\end{subequations} 
To close the bootstrap,  we show that we can in fact choose $K$
sufficiently large and $\epsilon$ sufficiently small so that the above
inequalities hold independently of $T$ with $K\ge$ replaced by $K\ge/2$. 
This implies that for small data, solutions can be extended for all $T > 0$. 

Applying the energy estimate of Lemma \ref{L:SymEnergy} to \eqref{E:SymZ} 
and summing over
$|\ga| \leq k$, we have
\begin{equation*}
	 M_1(t) \leq C \left(C_{\text{data}}\ge +\int_0^t C_{N_1} N_1(s)^2M_1(s) \,ds \right) 
	 	\cdot \exp{\left(C\int_0^t N_1(\gt)^2 \,d\gt \right)},
	\end{equation*}

where $C_{N_1}$ is a constant absorbing possible ``extra'' factors of $N_1(s)$ and only reflects the finiteness of $N_1(s)$.    Likewise, applying Proposition \ref{P:DivIeq} to \eqref{E:DivZ} and Proposition \ref{P:Linfty} to \eqref{E:BoxZ} gives us
	\begin{equation*}
	 M_2(t) \leq C \ge + \int_0^t C_{N_1}\,N_1(s)^2 M_1(s)\,ds
	\end{equation*}
and
	\begin{equation*}
	 N_2(t) \leq \frac{C}{(1+t)^{(n-1)/2}} \left( C_{\text{data}}\ge
		+ \int_0^t \frac{(N_1(s)+N_2(s))}{(1+s)^{(n-1)/2}}\left( M_1(s)+M_2(s)\right)^2 ds\right)
	\end{equation*}
where we have made use of the Cauchy-Schwartz inequality.  Under the assumed
bounds, we have that for some $C$, 
\begin{align*}	
	M_1(t) &\leq e^{C(K\ge)^2} \left(C_{\text{data}} \ge + C (K\ge)^2 K\ge\right) \leq \frac{K\ge}{2} ,
		\\
	M_2(t) &\leq C_{\text{data}} \ge + C (K\ge)^2 K\ge \leq \frac{K\ge}{2} ,
		\\
	N_2(t) & \leq \frac{1}{(1+t)^{(n-1)/2}} \left(C\,C_{\text{data}} \ge 
+ C (K\ge)^4 \right) 
		\leq \frac{K\ge}{2(1+t)^{(n-1)/2}},
\end{align*} 
where the second inequality in each line holds for all $t$ if  $K$ is chosen
sufficiently large and $\ge$ is chosen sufficiently small.   Recall that an
estimate for $N_1(t)$ follows from the estimate for $N_2(t)$.  Obtaining
these tighter bounds on the norms, we have closed the bootstrap. In view of
the continuation property stated in Theorem \ref{T:cauchy}, we have proved
small data global existence for $n \geq 3$.

\section{Global existence in dimension $n=2$}\label{S:DimTwo}
In the case of $n=2$, we require more detailed information concerning the
structure of the system.  In particular, we exploit the fact that the system
satisfies the so-called ``null-condition'' of Klainerman \cite{Kl_Null},
\cite{Kl_MoreNull}, which is a condition on the quadratic part of the
nonlinearity.  We consider null forms, quadradic forms of first derivatives,
which are given by
	\begin{equation}
	Q_{00}(f,g) = \gh^{\gm\gn} (\delmu f)(\delnu g), \qquad Q_{\ga\gb}(f,g) =  (\dela f)(\delb g)-(\delb f)(\dela g), \quad \ga \neq \gb
	\end{equation}
and satisfy time decay closer to that of cubic terms.  In particular, if $Q$ is any null form, then
	\begin{equation}\label{E:NullEst}
	\lvert  Q(f,g)(t,x) \rvert 
	\leq \frac{C}{1+t+|x|} \sum_{|\ga|=1}\lvert Z^{\ga}f (t,x)\rvert \sum_{|\gb|=1}\lvert Z^{\gb}g (t,x)\rvert. 
	\end{equation}
Furthermore, for any Lorentz vector field $Z$ and null form $Q$ there exists constants $a^{\gm\gn}$ so that
	\begin{equation}
	ZQ(f,g) = Q(Zf,g) + Q(f,Zg) + a^{\gm\gn}Q_{\gm\gn}(f,g).
	\end{equation}

Returning to the system \eqref{ELEqn}, we note that the Lagrangian associated to the volume element of the induced metric is $\scL = \sqrt{-\det{h}}$.  For small $|\del{}f|$, we have 
 \begin{align*}
	 -\det{h} &= 1 + \gh^{\gm\gn}\gd_{IJ}f^I_{\gm}f^J_{\gn} + \calO(|\del{}f|^4)\\
	 	&= 1 + \gd_{IJ} Q_{00}(f^I,f^J) + \calO(|\del{}f|^4)
	\end{align*}
and thus the Euler-Lagrange equations take the form
	\begin{equation}
	(1+ \gd_{KL}Q_{00}(f^K,f^L)) \Box f^J = \frac12\gh^{\gm\gn}f^J_{\gm} \delnu \left[ \gd_{AB}Q_{00}(f^A,f^B) \right] 
		+ \calO\left( |\Del{2}f|\cdot |\del{}f|^4 \right).
	\end{equation}
For small $|\del{}f|$ we have that
	\begin{equation}
	\left( 1+ \gd_{KL}Q_{00}(f^K,f^L)\right)^{-1} = 1 + \calO\left(|\del{}f|^2 \right) ;
	\end{equation}
thus we obtain
	\begin{equation}
	 \Box f^J  = \frac12 Q_{00}\left(f^J,\gd_{AB}Q_{00}(f^A,f^B)\right) + \calO\left( |\Del{2}f|\cdot |\del{}f|^4 \right).
	\end{equation}
%
Applying $Z^{\ga}$ yields
	\begin{equation}\label{E:NullZ}
	\Box (Z^{\ga} f^J) =  \sum_{\sum |\ga_i|\leq |\ga|} Q_{00}(Z^{\ga_1}f^{J},\gd_{AB}Q_{00}(Z^{\ga_2}f^{A},Z^{\ga_3}f^{B)})
		+ \calO\left( |Z^{\gb_1}\Del{2}f|\cdot |Z^{\gb_2}(\del{}f)^4| \right),
	\end{equation}	
where $ |\ga_1|+|\gb_2|\leq |\ga|$.

The proof of global existence when $n=2$ follows closely the proof for higher
dimensions, with three differences.  First, the null equation \eqref{E:NullZ} is
used in place of \eqref{E:BoxZ}.  Second, we use the following variation of
the energy estimate of Lemma \ref{L:SymEnergy}.
\begin{lemma}\label{L:EnIeqVar}
Under the hypotheses of Lemma \ref{L:SymEnergy} we have
	\begin{align*}
	\ltwot{t}{\del{}f} \leq &\,C \ltwot{0}{\del{}f} \cdot \exp{\left(-\int_0^t  C  \linftyt{\gt}{\del{}H}\,d\gt\right)}\\
	&+ C \int_0^t\sum_I\ltwot{s}{B_I}
		 \cdot\exp{\left(-\int_s^tC\linftyt{\gt}{\del{}H}\,d\gt\right)}\,ds.
	\end{align*} 
\end{lemma}
The third difference is that the bootstrap assumptions of \eqref{E:NormBd} are replaced by the following:
	\begin{subequations}\label{E:NormzTwo}
	\begin{align}
	M_1(t) & \leq K\ge (1+t)^{\gd},\\
	M_2(t) & \leq K\ge (1+t)^{\gd},\\
	N_1(t) & \leq \frac{K\ge}{(1+t)^{1/2}},\\
	N_2(t) & \leq \frac{K\ge}{(1+t)^{1/2}},
	\end{align} 
	\end{subequations}
where $0<\gd<\frac{1}{2}$ 
is a fixed, arbitrary constant.  We apply Lemma \ref{L:EnIeqVar} to \eqref{E:SymZ} and obtain the estimate
	\begin{equation*}
	 M_1(t) \leq  C_{\text{data}} \ge \cdot \exp{\left(\int_0^t \frac{C\, (K\ge)^2}{(1+\gt)}\,d\gt \right)}
	 	+ \int_0^t \frac{C\, (K\ge)^3}{(1+s)^{1-\gd}}
	 	\cdot \exp{\left(\int_s^t \frac{C\, (K\ge)^2}{(1+\gt)}\,d\gt \right)} ds.
	\end{equation*}
Computing
	\begin{equation*}
	 \exp{\left(\int_s^t \frac{C\, (K\ge)^2}{(1+\gt)}\,d\gt \right)}
	 	= \exp{\left(C\,(K\ge)^2 \log{\left(\frac{1+t}{1+s}\right)}\right)} = \left(\frac{1+t}{1+s}\right)^{C(K\ge)^2},
	\end{equation*}
one sees that
	\begin{equation*}
	\begin{aligned}
	 M_1(t) &\leq  C_{\text{data}} \ge (1+t)^{C(K\ge)^2} 
	 	+ C(K\ge)^3(1+t)^{C(K\ge)^2} \int_0^t (1+s)^{\gd-1-C(K\ge)^2}\,ds\\
	& = C_{\text{data}} \ge (1+t)^{C(K\ge)^2} + C(K\ge)^3(1+t)^{C(K\ge)^2}\left[(1+t)^{\gd-C(K\ge)^2}-1\right]\\
	&\leq K\ge (1+t)^{\gd}/2
	\end{aligned}
	\end{equation*}
for large $K$ and small $\ge$.  Similarly, application of Proposition 
\ref{P:DivIeq} to \eqref{E:DivZ} implies that for suitable large $K$ and small $\ge$ we have
that
	\begin{equation*}
	 \begin{aligned}
	M_2(t) &\leq C_{\text{data}} \ge \log{(2+t)} + \int_0^t \frac{C\, (K\ge)^3}{(1+s)^{1-\gd}}\,ds\\
		&\leq C_{\text{data}} \ge \log{(2+t)} + \tilde{C} (K\ge)^3(1+t)^{\gd}  \leq (1+t)^{\gd}{K\ge}/{2}.
	\end{aligned} 
	\end{equation*}
Lastly, via application of Proposition \ref{P:Linfty} to \eqref{E:NullZ} and use of the null estimate \eqref{E:NullEst} we obtain
	\begin{equation*}
	\begin{aligned}
	N_2(t) & \leq \frac{C}{(1+t)^{1/2}} \left( C_{\text{data}} \ge
		+\int_0^t \frac{C\,(K\ge)^3}{(1+s)^{3/2-2\gd}}\,ds \right) \\
		&\leq \frac{1}{(1+t)^{1/2}} \left(C\,C_{\text{data}} \ge + \check{C} (K\ge)^3 \right) \\
		&\leq \frac{K\ge}{2(1+t)^{1/2}},
	\end{aligned}
	\end{equation*}

provided $K$ is large and $\ge$ is small.  As the estimate for $N_2(t)$ implies the desired estimate for $N_1(t)$, we have that for  all finite intervals $[0,T)$, the estimates \eqref{E:NormzTwo} imply that better estimates hold, thus concluding the proof.\\ \\

\section{Concluding remarks}
The analysis of timelike minimal surfaces and submanifolds introduces a
family of geometrically-motivated quasilinear PDE systems which are both
intriguing and relatively unstudied. The work done by Brendle and Lindblad
and continued  here constitutes merely a first step in this analysis.

Among other things, one would hope to see progress on the characterization of
local well-posedness for  the timelike minimal submanifold PDEs in general
Lorentz spaces, as well as discovery of  nontrivial stable solutions other
than the flat planes one finds in Minkowski space. Exploring  the nature and
formation of singularities in these submanifolds should lead to interesting
new phenomena.

Further, it is interesting to explore the relation between the
timelike minimal surface PDE's and the classical minimal surface problem,
see e.g. \cite{TW:bernstein} and references therein.

In addition to the case considered in this paper, 
the more general $D$-brane equations which 
arise when various matter fields are incorporated into the action 
occur in string theory and are of current interest in theoretical physics.    
It is clear that timelike minimal submanifolds will provide a rich source of new mathematical problems for some time to come.
   
%



\subsection*{Acknowledgements}
We thank Mark Keel and Igor Rodnianski for useful discussions. We thank the
MFO at Oberwolfach and the Isaac Newton Institute in Cambridge for providing
stimulating places in which to collaborate on this research.  This work is
partially supported by NSF grants PHY-0354659 at Oregon, and DMS-0407732 at
the University of Miami.

	



\end{document}